\documentclass[11pt]{article}
\usepackage{epsfig}
\usepackage{amssymb,amsmath,amsthm,amscd}
\usepackage{latexsym}

\newtheorem{thm}{Theorem}[section]
\newtheorem{prop}[thm]{Proposition}

\newtheorem{lem}[thm]{Lemma}

\theoremstyle{definition}
\newtheorem{defn}[thm]{Definition}

\newcounter{labelflag} 

\newcommand{\Label}[1]{
                       \ifnum\thelabelflag=1
                          \ifmmode
                             \makebox[0in][l]{\qquad\fbox{\rm#1}}
                          \else
                             \marginpar{\vspace{0.7\baselineskip}
                                        \hspace{-1.1\textwidth}
                                        \fbox{\rm#1}}
                          \fi
                       \fi
                       \label{#1}
                      }
\newcommand{\be}{\begin{equation}}
\newcommand{\ee}{\end{equation}}

\newcommand{\h}{L^2(\R^n)}

\newcommand{\ii}{\int_{\R^n}}
\newcommand{\zto}{z(\theta_t \omega) }
\newcommand{\rh }{\rho ({\frac {|x|^2}{k^2}})}
\newcommand{\rhp }{\rho^\prime ({\frac {|x|^2}{k^2}})}

 \newcommand{\R}{\mathbb{R}}

\begin{document}

\begin{titlepage}
\title{\Large\bf   Upper Semicontinuity of
Random Attractors for  Non-compact Random Dynamical Systems}
\vspace{7mm}

\author{
Bixiang Wang  \thanks {Supported in part by NSF  grant DMS-0703521}
\vspace{1mm}\\
Department of Mathematics, New Mexico Institute of Mining and
Technology \vspace{1mm}\\ Socorro,  NM~87801, USA \vspace{1mm}\\
Email: bwang@nmt.edu}
\date{}
\end{titlepage}

\maketitle

\medskip

\begin{abstract}
The upper semicontinuity of random attractors
for non-compact random dynamical systems is proved
when the union of all perturbed random attractors
is precompact with probability one. This result is applied to the
stochastic Reaction-Diffusion  with white noise
defined on the entire space  $\R^n$.
\end{abstract}

{\bf Key words.}    Stochastic Reaction-Diffusion equation,
random attractor,
asymptotic compactness, upper semicontinuity.

 {\bf MSC 2000.} Primary 37L55. Secondary 60H15, 35B40.

%

%
%

\baselineskip=1.5\baselineskip

\section{Introduction}
\setcounter{equation}{0}
In this paper,  we  study the limiting behavior of random attractors
of non-compact random dynamical systems as stochastic perturbations
approach zero. In particular, we will establish the upper semicontinuity
of random attractors for the stochastically perturbed
Reaction-Diffusion equation defined on the entire space $\R^n$:
\be
\label{intro1}
du + (\lambda u - \Delta u )dt = (f(x,u) + g(x)) dt +  \epsilon h  dW,
 \ee
 where  $\epsilon$ is a small positive parameter,
 $\lambda$ is a fixed positive constant, $g$ and $h $
 are given functions
 defined on $\R^n$, $f$ is a smooth nonlinear  function satisfying
  some  conditions, and
 $W$ is a
 two-sided real-valued  Wiener process on a complete probability
 space.

 By a random attractor  we mean a compact  and invariant random set
 which attracts all solutions when initial times approach minus infinity.
 The concept of random attractor
 was   introduced in \cite{cra2, fla1} as
   extension to    stochastic  systems
 of  the  concept of global attractor   for
  deterministic   equations found in \cite{bab1, hal1, rob1, sel1, tem1}, for instance.
  In the case of bounded domains, random attractors
  for stochastic PDEs have been studied by many authors, see, e.g.,
  \cite{car1, car3, car4, car5, car6,  cra1, cra2, fla1, kap1, klo1, lv1, lv2, rob2,
  sch1, yan1, zho1} and the
references therein.
In these papers, the  asymptotic compactness of random dynamical systems
follows directly from the compactness of Sobolev embeddings  in bounded domains.
This is the  key  to prove the existence of random attractors
for PDEs defined in bounded domains.
Since Sobolev embeddings are not compact on unbounded domains,
the random dynamical systems associated with   PDEs
in this case are non-compact, and the asymptotic compactness of solutions
cannot be obtained simply   from  these embeddings.
This is a reason why there are only a few  results on existence of random attractors
for PDEs defined on unbounded domains. Nevertheless, the existence of
such attractors for some stochastic PDEs on unbounded domains has
been proved in \cite{blw, wan2} recently.
The   asymptotic compactness and existence of   absorbing sets
 for the
   stochastic Navier-Stokes equations
 on unbounded domains were established in    \cite{brz1}.

 In this paper, we will examine the limiting behavior of random attractors
 for the stochastically perturbed Reaction-Diffusion
 equation \eqref{intro1} defined on $\R^n$ when $\epsilon \to 0$, and prove
 the upper semicontinuity of these perturbed random attractors.
 In the deterministic case,
 the upper semicontinuity of global attractors
 were investigated in \cite{hal1, hal2, hal3, tem1}  and many others.
 For stochastic PDEs defined in bounded domains, this problem has been
 studied by the authors of \cite{car5,car6, lv1, lv2, rob2}.
 To the best of our knowledge, there is no result reported
 in the literature on the upper semicontinuity
 of random attractors for stochastic PDEs defined on unbounded domains.
 The purpose of this paper is to prove such a result for equation \eqref{intro1}
 on $\R^n$. Of course,  the main difficulty here is the non-compactness of
 Sobolev embeddings on $\R^n$.
 In this paper, we will overcome the obstacles caused by
 the non-compactness of embeddings by using  uniform estimates
 for far-field values of functions inside the perturbed random attractors.
 Actually, by a cut-off technique, we will show that
 the values of all functions in all perturbed random attractors
 are uniformly convergent to zero (in a sense) when spatial variables
 approach infinity (see the proof of  Lemma \ref{lem61} for more details).

 The outline of this paper is  as follows.
 We recall the basic random attractors theory in the
 next section, and prove a result on the upper semicontinuity
 of random attractors in Section 3. This result works for
 non-compact random dynamical systems corresponding
 to stochastic PDEs defined on unbounded domains.
 In Section 4, we define a continuous random dynamical system
 for equation \eqref{intro1} in $L^2(\R^n)$. The uniform estimates
 of solutions for the equation are given in Section 5.
 Finally, we prove the upper semicontinuity of random attractors
 for \eqref{intro1}  in the last section.

  We denote by
$\| \cdot \|$ and $(\cdot, \cdot)$ the norm and the inner product
in  $L^2(\R^n)$ and use $\| \cdot\|_{p}$    to denote   the norm  in
$L^{p}(\R^n)$.    Otherwise, the
norm of  a general  Banach space $X$  is written as    $\|\cdot\|_{X}$.
 The letters $c$ and $c_i$ ($i=1, 2, \ldots$)
are  generic positive constants  which may change their  values from line to
line or even in the same line.

\section{Random attractors}
\setcounter{equation}{0}

We recall some basic concepts
related to random attractors for stochastic dynamical
systems. The reader is referred to \cite{arn1, bat1, cra1, fla1} for more details.

Let  $(X, \| \cdot \|_X)$ be a
Banach  space with Borel $\sigma$-algebra $\mathcal{B}(X)$,
 and let
$(\Omega, \mathcal{F}, P)$  be  a probability space.

\begin{defn}
$(\Omega, \mathcal{F}, P,  (\theta_t)_{t\in \R})$
is called a metric  dynamical  system
if $\theta: \R \times \ \Omega \to \Omega$ is
$(\mathcal{B}(\R) \times \mathcal{F}, \mathcal{F})$-measurable,
$\theta_0$ is the identity on $\Omega$,
$\theta_{s+t} = \theta_t \circ  \theta_s$ for all
$s, t \in \R$ and $\theta_t P = P$ for all $t \in \R$.
 \end{defn}

\begin{defn}
\label{RDS}
A continuous random dynamical system (RDS)
on $ X  $ over  a  metric  dynamical system
$(\Omega, \mathcal{F}, P,  (\theta_t)_{t\in \R})$
is  a mapping
$$
\phi: \R^+ \times \Omega \times X \to X, \quad (t, \omega, x) \mapsto \phi(t, \omega, x),
$$
which is $(\mathcal{B}(\R^+) \times \mathcal{F} \times \mathcal{B}(X), \mathcal{B}(X))$-measurable and
satisfies, for $P$-a.e.  $\omega \in \Omega$,

(i) \  $\phi(0, \omega, \cdot) $ is the identity on $X$;

(ii) \  $\phi(t+s, \omega, \cdot) = \phi(t, \theta_s \omega, \cdot) \circ \phi(s, \omega, \cdot)$
for all $t, s \in \R^+$;

(iii) \  $\phi(t, \omega, \cdot): X \to  X$ is continuous for all $t \in  \R^+$.
\end{defn}

Hereafter, we always assume that
$\phi$ is a continuous RDS on $X$ over
$(\Omega, \mathcal{F}, P,  (\theta_t)_{t\in \R})$.

\begin{defn}
 A random  bounded set $\{B(\omega)\}_{\omega \in \Omega}$
 of  $  X$  is called  tempered
 with respect to $(\theta_t)_{t\in \R}$ if for $P$-a.e. $\omega \in \Omega$,
 $$ \lim_{t \to \infty} e^{- \beta t} \|  B(\theta_{-t} \omega ) \|_X =0
 \quad \mbox{for all} \  \beta>0,
 $$
 where $\| B\|_X =\sup\limits_{x \in B} \| x \|_{X}$.
\end{defn}

\begin{defn}
Let $\mathcal{D}$ be a collection of  random  subsets of $X$.
Then  $\mathcal{D}$ is called inclusion-closed if
   $D=\{D(\omega)\}_{\omega \in \Omega} \in {\mathcal{D}}$
and  $\tilde{D}=\{\tilde{D}(\omega)\}_{\omega \in \Omega} $
with
  $\tilde{D}(\omega) \subseteq D(\omega)$ for all $\omega \in \Omega$ imply
  that  $\tilde{D} \in {\mathcal{D}}$.
  \end{defn}

\begin{defn}
Let $\mathcal{D}$ be a collection of random subsets of $X$ and
$\{K(\omega)\}_{\omega \in \Omega} \in \mathcal{D}$. Then
$\{K(\omega)\}_{\omega \in \Omega} $ is called a random absorbing
set for   $\phi$ in $\mathcal{D}$ if for every $B \in \mathcal{D}$
and $P$-a.e. $\omega \in \Omega$, there exists $T(B, \omega)>0$ such
that
$$
\phi(t, \theta_{-t} \omega, B(\theta_{-t} \omega)) \subseteq K(\omega)
\quad \mbox{for all} \ t \ge T(B, \omega).
$$
\end{defn}

\begin{defn}
Let $\mathcal{D}$ be a collection of random subsets of $X$. Then
$\phi$ is said to be  $\mathcal{D}$-pullback asymptotically
compact in $X$ if  for $P$-a.e. $\omega \in \Omega$,
$\{\phi(t_n, \theta_{-t_n} \omega,
x_n)\}_{n=1}^\infty$ has a convergent  subsequence  in $X$
whenever
  $t_n \to \infty$, and $ x_n\in   B(\theta_{-t_n}\omega)$   with
$\{B(\omega)\}_{\omega \in \Omega} \in \mathcal{D}$.
\end{defn}

\begin{defn}
Let $\mathcal{D}$ be a collection of random subsets of $X$.
Then a random set $\{\mathcal{A}(\omega)\}_{\omega \in \Omega}$
of $X$
is called a   $\mathcal{D}$-random   attractor
(or $\mathcal{D}$-pullback attractor)  for
  $\phi$
if the following  conditions are satisfied, for $P$-a.e. $\omega \in \Omega$,

(i) \  $\mathcal{A}(\omega)$ is compact,  and
$\omega \mapsto d(x, \mathcal{A}(\omega))$ is measurable for every
$x \in X$;

(ii) \ $\{\mathcal{A}(\omega)\}_{\omega \in \Omega}$ is invariant, that is,
$$ \phi(t, \omega, \mathcal{A}(\omega)  )
= \mathcal{A}(\theta_t \omega), \ \  \forall \   t \ge 0;
$$

(iii) \ \ $\{\mathcal{A}(\omega)\}_{\omega \in \Omega}$
attracts  every  set  in $\mathcal{D}$,  that is, for every
 $B = \{B(\omega)\}_{\omega \in \Omega} \in \mathcal{D}$,
$$ \lim_{t \to  \infty} dist (\phi(t, \theta_{-t}\omega, B(\theta_{-t}\omega)), \mathcal{A}(\omega))=0,
$$
where $dist(\cdot, \cdot)$ is the Hausdorff semi-metric given by
$dist(Y,Z) =
  \sup\limits_{y \in Y }
\inf\limits_{z\in  Z}  \| y-z\|_{X}
 $ for any $Y\subseteq X$ and $Z \subseteq X$.
\end{defn}

The following existence result for a    random attractor
for a  continuous  RDS
can be found in \cite{bat1,blw,  fla1}.
\begin{prop}
\label{att} Let $\mathcal{D}$ be an inclusion-closed  collection of random subsets of
$X$ and $\phi$ a continuous RDS on $X$ over $(\Omega, \mathcal{F},
P,  (\theta_t)_{t\in \R})$. Suppose  that $\{K(\omega)\}_{\omega
\in K} $ is a closed  random absorbing set for  $\phi$  in
$\mathcal{D}$ and $\phi$ is $\mathcal{D}$-pullback asymptotically
compact in $X$. Then $\phi$ has a unique $\mathcal{D}$-random
attractor $\{\mathcal{A}(\omega)\}_{\omega \in \Omega}$ which is
given by
$$\mathcal{A}(\omega) =  \bigcap_{\tau \ge 0} \  \overline{ \bigcup_{t \ge \tau} \phi(t, \theta_{-t} \omega, K(\theta_{-t} \omega)) }.
$$
\end{prop}

\section{Upper semicontinuity of random attractors}
\setcounter{equation}{0}

In this section, we establish the upper semicontinuity of random
attractors when small random perturbations approach zero.
Let $(X, \|\cdot \|_X)$ be a Banach space and  $\Phi$ be an autonomous dynamical
system defined on $X$. Given $\epsilon>0$, suppose $\Phi_\epsilon$
is a random dynamical system over a metric system
$(\Omega, {\mathcal{F}}, P, (\theta_t)_{t\in \R} )$.
We further suppose that for $P$-a.e. $\omega \in \Omega$, $t\ge 0$,
$\epsilon_n \to 0$, and
$x_n$, $x \in X$ with $x_n \to x$, the following holds:
\be
\label{upp1}
\lim_{n \to \infty}
\Phi_{\epsilon_n} (t, \omega,  x_n )
= \Phi(t) x.
\ee

Let ${\mathcal{D}}$ be a collection of subsets of $X$. Given
$\epsilon>0$, suppose that $\Phi_\epsilon$ has a  random attractor
${\mathcal{A}}_\epsilon =\{\mathcal{A}_\epsilon (\omega)\}_{\omega \in \Omega}
\in {\mathcal{D}}$
and a random absorbing
set $E_\epsilon = \{E_\epsilon (\omega)\}_{\omega \in \Omega}
\in {\mathcal{D}}$  such that for some deterministic
positive constant $c$ and  for $P$-a.e. $\omega \in \Omega$,
\be
\label{upp2}
\limsup_{\epsilon \to  0}
\| E_\epsilon (\omega) \|_X \le c,
\ee
where $\| E_\epsilon (\omega) \|_X =
\sup\limits_{x \in  E_\epsilon (\omega)} \| x \|_X$.
We also assume that there exists $\epsilon_0 >0$ such that
for $P$-a.e. $\omega \in \Omega$,
\be
\label{upp3}
\bigcup_{0<\epsilon \le \epsilon_0} {\mathcal{A}}_\epsilon (\omega)
\quad \mbox{is precompact in } \quad X.
\ee
Let ${\mathcal{A}}_0$ be the global attractor of $\Phi$ in $X$, which means
that ${\mathcal{A}}_0$ is compact and invariant and attracts every bounded
subset of $X$ uniformly.
 Then the
relationships between ${\mathcal{A}}_\epsilon$ and ${\mathcal{A}}_0$
are given by the following theorem.

\begin{thm}
\label{semicont}
Suppose \eqref{upp1}-\eqref{upp3} hold. Then for $P$-a.e. $\omega \in \Omega$,
\be
\label{upp4}
dist ({\mathcal{A}}_\epsilon (\omega), {\mathcal{A}}_0 )
\to 0,
\quad \mbox{as} \quad \epsilon \to  0.
\ee
\end{thm}

\begin{proof}
We argue by contradiction. If \eqref{upp4} is not true, then there
$\delta>0$ and a sequence $\{x_n\}_{n=1}^\infty$
 with $x_n \in {\mathcal{A}}_{\epsilon_n} (\omega)$
and $\epsilon_n \to 0$ such that
\be
\label{upp55}
dist(x_n, \mathcal{A}_0) \ge \delta.
\ee
It follows from \eqref{upp3} that there are $y_0 \in X$
and a subsequence of $\{x_n\}_{n=1}^\infty$
(still denoted by $\{x_n\}_{n=1}^\infty$) such that
\be
\label{upp5}
\lim_{n \to \infty} x_n =y_0.
\ee
Next we prove $y_0 \in {\mathcal{A}}_0$. To this end, we
 take a sequence $\{t_m\}_{m=1}^\infty$ with $t_m \to \infty$.
By the invariance of ${\mathcal{A}}_{\epsilon_n}$ we find that
there exists a sequence
$\{x_{1,n}\}_{n=1}^\infty$ with $x_{1,n} \in
 {\mathcal{A}}_{\epsilon_n}(\theta_{-t_1} \omega)$
such that
\be
\label{upp6}
x_n = \Phi_{\epsilon_n} (t_1, \theta_{-t_1} \omega, x_{1,n} )  ,
\ \forall  \ n \ge 1.
\ee
By \eqref{upp3} again, there exist
$y_1 \in X$ and a subsequence of $\{x_{1,n}\}_{n=1}^\infty$
(still denoted by $\{x_{1,n}\}_{n=1}^\infty$) such that
\be
\label{upp7}
\lim_{n \to \infty} x_{1,n} = y_1.
\ee
By \eqref{upp1} and \eqref{upp7} we find that
\be
\label{upp8}
\lim_{n \to \infty}
\Phi_{\epsilon_n} (t_1, \theta_{-t_1} \omega, x_{1,n})
= \Phi (t_1) y_1.
\ee
It follows from \eqref{upp5}-\eqref{upp6}
and \eqref{upp8} that
\be
\label{upp9}
y_0 = \Phi(t_1) y_1.
\ee
Since
$ x_{1,n}  \in  {\mathcal{A}}_{\epsilon_n}(\theta_{-t_1} \omega)$
and
$ {\mathcal{A}}_{\epsilon_n}(\theta_{-t_1} \omega)
\subseteq { {E}}_{\epsilon_n}(\theta_{-t_1} \omega)$,
by \eqref{upp2} we get
\be
\label{upp10}
\limsup_{n \to \infty} \| x_{1,n} \|_X
\le \limsup_{n \to \infty}
\| E_{\epsilon_n} (\theta_{-t_1} \omega) \|_X
\le c.
\ee
By \eqref{upp7} and
\eqref{upp10} we find that
$$
 \| y_1 \|_X \le c.
$$
 Similarly, for each $m \ge 2$, repeating the above procedure, we can find that
 there is $y_m \in X$ such that
 \be
 \label{upp12}
 y_0 = \Phi(t_m) y_m, \quad \forall \ m \ge 2,
 \ee
 and
 \be
 \label{upp13}
 \| y_m \|_X \le c, \quad \forall \ m \ge 2.
 \ee
 Since $t_m \to \infty$, \eqref{upp12} and
 \eqref{upp13} imply that
 $y_0 \in {\mathcal{A}}_0$.
 Therefore, by \eqref{upp5} we have
 $$
 dist(x_n, {\mathcal{A}})
 \le dist(x_n, y_0) \to 0,$$
 a contradiction with \eqref{upp55}.
 This completes the proof.
\end{proof}

We remark that the upper semicontinuity of random attractors
for stochastic PDEs as perturbations of autonomous, non-autonomous
and random systems
was first proved by the authors in
\cite{car5}, \cite{car6} and \cite{rob2}, respectively.
The conditions
  \eqref{upp1}-\eqref{upp3} of this paper
are close  but different from that given in \cite{car5,car6, rob2}.
For instance, the following condition is  essentially assumed in \cite{car5, car6, rob2}
(see  Theorem 2 on page 1562 in \cite{car5}, Theorem 3.1 on page 496 in \cite{car6},
and Theorem 2 on page 655 in \cite{rob2}):
there exists a compact set $K$ such that, $P$-a.s.
\begin{equation}
\label{add_upp1}
\lim_{\epsilon \to 0}
dist ({\mathcal{A}}_\epsilon (\omega), K )
=0 .\end{equation}
For parabolic PDEs defined in {\it bounded} domains, the solution operators are compact,
which follows from the regularity of solutions and the compactness of
  Sobolev embeddings.
  In that case, the existence of the compact set $K$ satisfying
  condition \eqref{add_upp1} can be obtained by the existence of bounded absorbing sets
  in a space with higher regularity (see \cite{car5, car6, rob2}).
  However, this method does not work for PDEs defined on {\it unbounded} domains because
  Sobolev embeddings are no longer compact in this case.
  Therefore, in the case of unbounded domains, it is difficult to find a compact set
  $K$  which satisfies  \eqref{add_upp1}.
  In this paper, we require condition \eqref{upp3} rather than \eqref{add_upp1}.
  As proved in Section 6 of this paper, the condition \eqref{upp3}  is indeed fulfilled
  for the parabolic equation \eqref{intro1} defined on  the unbounded domain
  $\R^n$, and
  hence the upper semicontinuity of the random attractors
   follows from Theorem \ref{semicont} immediately.

\section{Stochastic Reaction-Diffusion equations
on  $\R^n$  }
\setcounter{equation}{0}

In this paper, we will investigate the upper semicontinuity
of   random attractors
of the  stochastic Reaction-Diffusion equation    defined on
$\R^n$.
Given a small positive parameter $\epsilon$, consider the
following stochastically perturbed equation:
 \be
\label{rd1}
du + (\lambda u - \Delta u )dt = (f(x,u) + g(x)) dt +  \epsilon  h\; dW, \quad  x \in \R^n, \ \ t >0,
 \ee
 with the initial condition:
 \be
 \label{rd2}
 u(x, 0) =u_0(x), \quad x \in \R^n.
 \ee
 Here $\epsilon$ and  $\lambda$ are   positive constants, $g $ is a given function  in
$L^2(\R^n)$,    $h  \in H^2(\R^n) \cap W^{2, p}(\R^n)$
for some  $p\ge 2$,     $W$   is a  two-sided real-valued
 Wiener process
 on a  complete probability
 space  $(\Omega,
\mathcal{F}, P)$, where $P$ is the Wiener distribution,
$ \Omega $ is a subset of
$   \{ \omega   \in C(\R, \R ): \
\omega(0) =  0 \}
$ with $P(\Omega)=1$,  and
$\mathcal{F}$ is a  $\sigma$-algebra.
In addition, the space $(\Omega,
\mathcal{F}, P)$ is invariant under the Wiener shift:
$$ \theta_t \omega (\cdot) = \omega (\cdot +t) - \omega (t), \quad  \omega \in \Omega, \ \ t \in \R .
$$
This means that $(\Omega, \mathcal{F}, P, (\theta_t)_{t\in \R})$ is a
metric  dynamical  system (see, e.g.,  \cite{car5, sch1} for existence of this space).

Consider  the one-dimensional Ornstein-Uhlenbeck
equation:
\be
\label{z1}
dy  + \lambda y  dt = dW(t).
\ee
One may easily check that a solution to \eqref{z1} is given by
$$
y (\theta_t \omega )
 = -\lambda \int^0_{-\infty} e^{\lambda \tau}  ( \theta_t \omega  ) (\tau) d \tau, \quad t \in \R.
$$
Note that the random variable $|y (\omega )|$ is tempered and $y (\theta_t\omega )$
is $P$-a.e. continuous. Therefore, it follows from Proposition 4.3.3 in \cite{arn1}
that there exists a tempered function $r(\omega)>0$ such that
\be
\label{z2}
      |y (\omega )|^2 + |y (\omega )|^p
\le r(\omega),
\ee
where $r(\omega)$ satisfies, for $P$-a.e. $\omega \in \Omega$,
\be
\label{z3}
r(\theta_t \omega ) \le e^{{\frac \lambda{2}} |t|} r(\omega), \quad t\in \R.
\ee
Then it follows from \eqref{z2}-\eqref{z3} that, for $P$-a.e. $\omega \in \Omega$,
\be
\label{zz}
   |y (\theta_t \omega )|^2 + |y (\theta_t \omega )|^p
\le e^{{\frac \lambda{2}} |t|} r(\omega), \quad t\in \R.
\ee
Let  $z(\theta_t \omega) =   h  y(\theta_t \omega  ) $
 and $v(t) = u(t) -  \epsilon z(\theta_t \omega)$ where $u$ is  a solution of  problem \eqref{rd1}-\eqref{rd2}.
Then $v$ satisfies
\be
\label{v1}
{\frac {\partial v}{\partial t}} + \lambda v -\Delta v = f(x, v + \epsilon z(\theta_t \omega) ) + g  + \epsilon  \Delta z (\theta_t \omega ).
\ee

In this paper, we assume that the nonlinearity
 $f$ satisfies  the following conditions:
 For all $x \in  R^n$ and $ s \in R$,
\begin{equation}
\label{f1}
f(x, s) s \le -  \alpha_1 |s|^p + \psi_1(x),
\end{equation}
\begin{equation}
\label{f2}
|f(x, s) |   \le \alpha_2 |s|^{p-1} + \psi_2 (x),
\end{equation}
\begin{equation}
\label{f3}
{\frac {\partial f}{\partial s}} (x, s)   \le \beta,
\end{equation}
\begin{equation}
\label{f4}
| {\frac {\partial f}{\partial x}} (x, s) | \le  \psi_3(x),
\end{equation}
where $\alpha_1$, $\alpha_2$ and    $\beta$
 are positive constants,
$\psi_1 \in L^1(R^n) \cap L^\infty(R^n)$, and $\psi_2 \in L^2(R^n) \cap L^q(\R^n) $
with ${\frac 1q} + {\frac 1p} =1$, and $\psi_3 \in L^2(\R^n)$.

It follows from \cite{blw} that, under conditions \eqref{f1}-\eqref{f4},
for $P$-a.e. $\omega \in \Omega $ and for all $v_0  \in L^2(\R^n)$,
   \eqref{v1}
has   a unique solution
$v(\cdot, \omega, v_0) \in C([0, \infty), L^2(\R^n)) \bigcap L^2((0,T), H^1(\R^n))$
  with $v(0, \omega, v_0) = v_0$ for every $T>0$.
Furthermore , the solution is continuous with respect to $v_0$ in
 $L^2(\R^n)$ for all $t \ge 0$.
 Let
 \be
 \label{uandv}
 u(t, \omega, u_0)= v(t,
\omega,   v_0) + \epsilon z(\theta_t \omega),
\quad \mbox{where} \quad v_0 = u_0 - \epsilon z(\omega).
\ee
We  can  associate a random dynamical system $\Phi_\epsilon$
with problem \eqref{rd1}-\eqref{rd2} via $u$ for each $\epsilon>0$, where
   $\Phi_\epsilon: \R^+ \times \Omega \times L^2(\R^n) \to
L^2(\R^n)$ is given by
 \be \label{phi}
 \Phi_\epsilon (t, \omega, u_0) =  u (t, \omega, u_0) , \quad \hbox{for every }\,\, \ (t, \omega, u_0) \in
\R^+ \times \Omega \times L^2(\R^n)  .
\ee
Then
$\Phi_\epsilon$
 is a continuous  random dynamical system
 over $(\Omega, \mathcal{F}, P, (\theta_t)_{t\in \R})$
 in $L^2(\R^n)  $.
 In the sequel, we always assume that
   $\mathcal{D}$  is  a collection of random subsets of
   $L^2(\R^n)$ given by
   \be
   \label{ddd}
   \mathcal{D} = \{ D=\{D(\omega)\}_{\omega \in \Omega}, \
D(\omega) \subseteq L^2(\R^n) \quad
\mbox{and} \quad e^{-{\frac 12} \lambda t} \| B(\theta_{-t} \omega ) \|
\to 0 \ \mbox{as} \  t \to \infty \},
\ee
where
 $$\| B(\theta_{-t} \omega) \|
 =\sup_{u \in  B(\theta_{-t} \omega) } \| u\|.
 $$
 In \cite{blw}, the authors proved that $\Phi_\epsilon$ has a $\mathcal{D}$-pullback
 random attractor if $\mathcal{D}$  is the collection of all tempered random subsets
 of $L^2(\R^n)$. Following the arguments of \cite{blw}, we can also
 prove that $\Phi_\epsilon$ has a unique $\mathcal{D}$-pullback
 random attractor $\{{\mathcal{A} }_\epsilon (\omega)\}_{\omega \in \Omega}$
  when $\mathcal{D}$  is given by  \eqref{ddd}
 (the existence of $\{{\mathcal{A}}_\epsilon (\omega)\}_{\omega \in \Omega}$
 in this case
 is also implied by the estimates given in Section 5 of this paper).
  When $\epsilon =0$, problem \eqref{rd1}-\eqref{rd2}
  defines  a continuous deterministic dynamical
  system $\Phi $  in $L^2(\R^n)$. In this case, the results  of \cite{blw}
  imply that
  $\Phi $ has a unique global attractor
     ${\mathcal{A}} $ in $L^2(\R^n)  $.
  The purpose of this paper is to establish  the relationships of
  $\{{\mathcal{A}}_\epsilon (\omega)\}_{\omega \in \Omega}
  $ and ${\mathcal{A}}$ when $\epsilon \to 0$.

\section{Uniform  estimates of solutions  }
\setcounter{equation}{0}

In this section, we
 derive uniform estimates of  solutions
 with respect to the small parameter  $\epsilon$. These estimates
  are
 useful
 for proving the semicontinuity of the perturbed random attractors.
Here and after,  we always assume that $\mathcal{D}$ is
 the collection of random  subsets of  $L^2(\R^n)  $
 given in \eqref{ddd}.

\begin{lem}
\label{lem41}  Let $0< \epsilon \le 1  $, $g \in L^2(\R^n)$ and
\eqref{f1}-\eqref{f4} hold.  Then for every
$B= \{ B(\omega) \}_{\omega \in \Omega}
\in {\mathcal{D}}$ and
   $P$-a.e. $\omega \in \Omega$, there
is $T(B,\omega)>0$,  independent of $\epsilon$, such that
for all
$v_0(\theta_{-t} \omega) \in B(\theta_{-t} \omega)$,
$$\|v(t, \theta_{-t} \omega, v_0(\theta_{-t} \omega) )\|^2
\le
e^{-\lambda t} \| v_0 ( \theta_{-t} \omega) \|^2
+    c + \epsilon c r(\omega), \quad \forall \ t \ge 0,
$$
$$
\int_0^t e^{\lambda (\tau -t)} \| \nabla v(\tau, \theta_{-t}\omega, v_0(\theta_{-t}\omega))\|^2 d\tau
\le
e^{-\lambda t} \| v_0 ( \theta_{-t} \omega) \|^2
+    c + \epsilon c r(\omega), \quad \forall \ t \ge 0,
$$
and
$$
 \int_0^t e^{\lambda (\tau -t)} \| u(\tau, \theta_{-t}\omega, u_0(\theta_{-t}\omega))\|^p_p d\tau \le  c + \epsilon c r(\omega),
 \quad \forall \ t \ge T(B, \omega),
$$
where $c$ is a positive deterministic constant independent
of $\epsilon$, and $r(\omega)$ is the tempered  function
  in \eqref{z2}.
\end{lem}

\begin{proof}
The  idea of proof  is similar to that given in \cite{blw}, but
now we have to pay attention to how the estimates depend on the
parameter $\epsilon$.
Multiplying \eqref{v1} by  $v$ and then integrating over $\R^n$, we find that
\be
\label{p41_1}
{\frac 12} {\frac d{dt}} \|v\|^2 + \lambda \| v\|^2 + \| \nabla v \|^2 =
\int_{\R^n} f(x, v+ \epsilon z(\theta_t \omega ) )  \ v  dx  + (g, v) +
 \epsilon (  \Delta z(\theta_t \omega ), v).
 \ee
For the nonlinear term, by \eqref{f1}-\eqref{f2} we obtain
$$
\int_{\R^n}  f(x, v+ \epsilon z(\theta_t \omega ) )  \  v  dx
=\int_{\R^n}  f(x, v+ \epsilon z(\theta_t \omega ) )   \
 ( v + \epsilon z(\theta_t \omega )  ) dx
- \epsilon \int_{\R^n}  f(x, v+ \epsilon z(\theta_t \omega ) ) \  \zto dx
$$
$$
\le -\alpha_1 \ii | u |^p dx + \ii \psi_1(x) dx
- \epsilon \ii f(x, u)\;  \zto dx
$$
\be
\label{p41_2}
\le - {\frac 12} \alpha_1 \| u \|^p_p
+ \epsilon c_2 ( \|\zto \|^p _p  +  \|\zto \|^2 ) +c_3,
\ee
where $c_2$ and $c_3$ do not depend on $\epsilon$.
Similarly, the remaining   terms on the
right-hand side of \eqref{p41_1} are bounded by
\be
\label{p41_3}
\| g \| \| v \| + \epsilon \|\nabla \zto \| \| \nabla v \|
\le
{\frac 12} \lambda \| v \|^2 + {\frac 1{2\lambda}} \| g \|^2
+{\frac 12}\epsilon  \| \nabla \zto \|^2
+{\frac 12}   \| \nabla v \|^2.
\ee
Then it follows from \eqref{p41_1}-\eqref{p41_3} that
\be
\label{p41_4}
 {\frac d{dt}} \|v\|^2 + \lambda \| v\|^2 + \| \nabla v \|^2
 + \alpha_1 \| u \|^p_p
 \le \epsilon  c_4 ( \|\zto \|^p _p  +  \|\zto \|^2 + \| \nabla \zto \|^2  ) +c_5.
\ee
Note that $\zto =   h   y(\theta_t \omega) $ and $h  \in H^2(\R^n) \cap W^{2,p}(\R^n)$.
Then we have
\be
\label{p41_5}
  \|\zto \|^p _p  +  \|\zto \|^2 + \| \nabla \zto \|^2
  \le
c_6 (  | y (\theta_t \omega )   |^p +  | y (\theta_t \omega )   |^2 )
=p_1(\theta_t \omega)  .
\ee
 By  \eqref{zz}, we find that for $P$-a.e. $\omega \in \Omega$,
 \be
 \label{p41_6a1}
 p_1(\theta_\tau \omega)
\le  c_6  e^{{\frac 12} \lambda |\tau|} r(\omega) , \quad \forall \ \tau \in \R.
\ee
It follows from \eqref{p41_4}-\eqref{p41_5} that,  for all $t \ge 0$,
\be
\label{p41_6}
 {\frac d{dt}} \|v\|^2 + \lambda \| v\|^2 + \| \nabla v \|^2
 + \alpha_1 \| u \|^p_p
 \le \epsilon c_4  p_1(\theta_t \omega)  + c_5.
\ee
Multiplying \eqref{p41_6} by $e^{\lambda t}$ and then integrating the
inequality, we get that, for all $t\ge 0$,
$$
\| v(t, \omega, v_0(\omega) ) \|^2
+ \int_0^t e^{\lambda (\tau -t)} \| \nabla v(\tau, \omega, v_0(\omega))\|^2 d\tau
+ \alpha_1 \int_0^t e^{\lambda (\tau -t)} \| u(\tau, \omega, u_0(\omega))\|^p_p d\tau
$$
\be
\label{p41_8}
\le e^{-\lambda t} \| v_0 (\omega) \|^2
+ \epsilon c_4 \int_0^t e^{\lambda(\tau -t)}  p_1(\theta_\tau \omega ) d\tau
+ c_7.
\ee
By replacing $\omega$ by $\theta_{-t} \omega$, we get from
\eqref{p41_8} and \eqref{p41_6a1}  that, for all $t \ge 0$,
$$
\| v(t,  \theta_{-t} \omega, v_0( \theta_{-t} \omega) ) \|^2
+ \int_0^t e^{\lambda (\tau -t)} \| \nabla v(\tau, \theta_{-t}\omega, v_0(\theta_{-t}\omega))\|^2 d\tau
+ \alpha_1 \int_0^t e^{\lambda (\tau -t)}
\| u(\tau, \theta_{-t}\omega, u_0(\theta_{-t}\omega))\|^p_p d\tau
$$
    \be
\label{p41_9}
\le e^{-\lambda t} \| v_0 ( \theta_{-t} \omega) \|^2
+ \epsilon c_4 \int_0^t e^{\lambda(\tau-t)}  p_1(  \theta_{\tau-t} \omega ) ds
+ c_7
\le
e^{-\lambda t} \| v_0 ( \theta_{-t} \omega) \|^2
+  \epsilon c_9   r(\omega)
+   {c_7} .
\ee
Since $v_0(\theta_{-t} \omega ) \in \mathcal{D}$, there is
$T=T(B, \omega)$, independent of $\epsilon$, such that for
all $t \ge T$,
$$e^{-\lambda t} \| v_0 ( \theta_{-t} \omega) \|^2
\le 1,
$$
 which along with
\eqref{p41_9} implies the lemma.
\end{proof}

As a consequence of Lemma \ref{lem41}, we have the following
estimates for $u$.

\begin{lem}
\label{lem42}
 Let $0< \epsilon \le 1  $, $g \in L^2(\R^n)$ and
\eqref{f1}-\eqref{f4} hold.  Then for every
$B= \{ B(\omega) \}_{\omega \in \Omega}
\in {\mathcal{D}}$ and
   $P$-a.e. $\omega \in \Omega$, there
is $T(B,\omega)>0$,  independent of $\epsilon$, such that
for all $t \ge T(B, \omega)$ and
$u_0(\theta_{-t} \omega) \in B(\theta_{-t} \omega)$,
$$\|u(t, \theta_{-t} \omega, u_0(\theta_{-t} \omega) )\|^2
\le  c + \epsilon c r(\omega),
$$
and
$$
\int_t^{t+1} \| \nabla u(\tau, \theta_{-t-1} \omega, u_0(\theta_{-t-1 } \omega)) \|^2 d\tau
\le c + \epsilon c  r  (\omega)  ,
$$
where $c$ is a positive deterministic constant independent
of $\epsilon$, and $r(\omega)$ is the tempered  function
  in \eqref{z2}.
\end{lem}

\begin{proof}
It follows from \eqref{uandv} and Lemma \ref{lem41} that
$$
\|  u(t, \theta_{-t} \omega, u_0(\theta_{-t} \omega) ) \|^2
\le 2 \| v(t,  \theta_{-t} \omega, u_0(\theta_{-t} \omega) - \epsilon z(\theta_{-t} \omega) ) \|^2
+ 2  \epsilon^2  \|  z (  \omega )  \|^2
$$
\be
\label{p42_20}
\le
4 e^{-\lambda t} ( \|  u_0(\theta_{-t} \omega) \|^2 +    \|  z(\theta_{-t} \omega)   \|^2 )
+ c + \epsilon  c  r(\omega),
\ee
where we have used \eqref{z2} and the fact $0<\epsilon \le 1$.
Since  $u_0(\theta_{-t} \omega ) \in B(\theta_{-t} \omega )$
and $\| z(\omega)\|^2$ is   tempered,   there is
 $T(B,\omega)>0$, independent of $\epsilon$,  such that for all $t \ge T(B, \omega)$,
 \be
 \label{p42_21}
  e^{-\lambda t} ( \|  u_0(\theta_{-t} \omega) \|^2 + \|  z(\theta_{-t} \omega)   \|^2 )
 \le 1,
 \ee
 which along with \eqref{p42_20} implies  that, for all $t \ge T (B,\omega)$,
\be
\label{p42_22}
\|  u(t, \theta_{-t} \omega, u_0(\theta_{-t} \omega) ) \|^2
\le
4 + c+  \epsilon c r(\omega)   .
 \ee
 Similarly, we have
 $$
 \| \nabla u  (\tau, \theta_{-t-1} \omega, u_0(\theta_{-t-1 } \omega)) \|^2
 =
  \| \nabla v (\tau, \theta_{-t-1} \omega, u_0(\theta_{-t-1 } \omega) -\epsilon z(\theta_{-t-1} \omega)  )
+ \epsilon \nabla z(\theta_{\tau-t-1} \omega) \|^2
 $$
 \be
 \label{p42a_1}
 \le
 2 \| \nabla v (\tau, \theta_{-t-1} \omega, u_0(\theta_{-t-1 }
  \omega) -\epsilon z(\theta_{-t-1} \omega)  ) \|^2
+ 2  \epsilon^2 \| \nabla z(\theta_{\tau-t-1} \omega) \|^2
  \ee
For $\tau \in (t, t+1)$, by \eqref{zz} we  find that
\be
\label{p42a_2}
  \|\nabla z(\theta_{\tau-t-1} \omega) \|^2
\le c |y (\theta_{\tau-t-1} \omega  ) |^2
\le c e^{{\frac \lambda 2}  } r (\omega).
\ee
By \eqref{p42a_1} and \eqref{p42a_2}, we get
$$
 \| \nabla u  (\tau, \theta_{-t-1} \omega, u_0(\theta_{-t-1 } \omega)) \|^2
\le
2 \| \nabla v (\tau, \theta_{-t-1} \omega, u_0(\theta_{-t-1 } \omega) -\epsilon z(\theta_{-t-1} \omega)  ) \|^2
+\epsilon c r(\omega).
$$
Integrating the above with respect to $\tau$
in $(t, t+1)$ we obtain
$$
 \int_t^{t+1} \| \nabla u  (\tau, \theta_{-t-1} \omega, u_0(\theta_{-t-1 } \omega)) \|^2 d \tau
 $$
 \be
 \label{p42a_3}
\le
2 \int_t^{t+1} \| \nabla v (\tau, \theta_{-t-1} \omega, u_0(\theta_{-t-1 } \omega) -\epsilon z(\theta_{-t-1} \omega)  ) \|^2 d\tau
+\epsilon c r(\omega).
\ee
Given $t \ge 0$, replacing  $t$ by $t+1$ in Lemma \ref{lem41} we find that
$$
\int_t^{t+1} e^{\lambda (\tau -t-1)} \|\nabla v (\tau, \theta_{-t-1} \omega, u_0(\theta_{-t-1 } \omega) -\epsilon z(\theta_{-t-1} \omega)  ) \|^2 d\tau
$$
\be
\label{p42a_4}
\le
2e^{-\lambda (t+1)} (\| u_0( \theta_{-t-1 } \omega )\|^2
 +    \| z(\theta_{-t-1} \omega )\|^2 )
+c +\epsilon c r(\omega).
\ee
Replacing $t$ by $t+1$ in \eqref{p42_21}, we find that
the first term on the right-hand side of \eqref{p42a_4}
is less than $2$ when $t \ge T(B,\omega)$. Therefore, we
have, for all $t \ge T(B, \omega)$,
$$
\int_t^{t+1} e^{\lambda (\tau -t-1)} \|\nabla v (\tau, \theta_{-t-1} \omega, u_0(\theta_{-t-1 } \omega) -\epsilon z(\theta_{-t-1} \omega)  ) \|^2 d\tau
\le 2 +c +\epsilon c r(\omega). $$
Since $ e^{\lambda (\tau -t-1)} \ge e^{-\lambda}$ for $\tau \in (t, t+1)$,
the above implies that,
for all $t \ge T(B, \omega)$,
\be
\label{p42a_5}
\int_t^{t+1}   \|\nabla v (\tau, \theta_{-t-1} \omega, u_0(\theta_{-t-1 } \omega) -\epsilon z(\theta_{-t-1} \omega)  ) \|^2 d\tau
\le e^{ \lambda} (2 +c +\epsilon c r(\omega)).
\ee
It follows from \eqref{p42a_3} and \eqref{p42a_5} that, for
all $t \ge T(B,\omega)$,
$$
 \int_t^{t+1} \| \nabla u  (\tau, \theta_{-t-1} \omega, u_0(\theta_{-t-1 } \omega)) \|^2 d \tau
 \le c + \epsilon c r(\omega),
$$
which along with \eqref{p42_22} concludes the proof.
\end{proof}

We are now in a position to establish the
uniform estimates of solutions in $H^1(\R^n)$.

\begin{lem}
\label{lem43}
Let $0< \epsilon \le 1  $, $g \in L^2(\R^n)$ and
\eqref{f1}-\eqref{f4} hold.  Then for every
$B= \{ B(\omega) \}_{\omega \in \Omega}
\in {\mathcal{D}}$ and
   $P$-a.e. $\omega \in \Omega$, there
is $T(B,\omega)>0$,  independent of $\epsilon$, such that
for all $t \ge T(B, \omega)$,
$u_0(\theta_{-t} \omega) \in B(\theta_{-t} \omega)$
and $v_0(\theta_{-t} \omega) = u_0(\theta_{-t} \omega) -\epsilon z(\omega)$,
$$\|\nabla v(t, \theta_{-t} \omega, v_0(\theta_{-t} \omega) )\|^2
\le  c + \epsilon c r(\omega),
$$
and
$$\|\nabla u(t, \theta_{-t} \omega, u_0(\theta_{-t} \omega) )\|^2
\le  c + \epsilon c r(\omega),
$$
where $c$ is a positive deterministic constant independent
of $\epsilon$, and $r(\omega)$ is the tempered  function
  in \eqref{z2}.
\end{lem}

\begin{proof}
Taking the inner product of  \eqref{v1}
  with $\Delta v$ in $L^2(\R^n)$, we get that
\be
\label{p43_1}
{\frac 12} {\frac d{dt}} \| \nabla v \|^2
+ \lambda \| \nabla  v \|^2
+ \| \Delta v \|^2
=-\ii f(x, u) \Delta v dx
-(g + \epsilon \Delta \zto, \Delta v ).
\ee
By \eqref{f2}-\eqref{f4}, the first term on the right-hand side
of   \eqref{p43_1} satisfies
$$
-\ii f(x, u) \;  \Delta v dx
= -\ii f(x, u) \;  \Delta u dx + \epsilon  \ii f(x, u) \; \Delta \zto dx
$$
$$
=\ii {\frac {\partial f}{\partial x}} (x, u) \; \nabla u dx
+ \ii {\frac {\partial f}{\partial u}} (x,u) \;  | \nabla  u |^2 dx
+ \epsilon \ii f(x, u) \Delta \zto dx
$$
\be
\label{p43_2}
\le c \left (
 \| \nabla u \|^2 + \| u \|^p_p \right )
 + \epsilon c \left (
  \|   \Delta \zto  \|^2 +  \|   \Delta \zto  \|^p_p
\right ) + c,
\ee
where we have used the fact  $0< \epsilon \le 1$.
For  the last term on the right-hand side
of \eqref{p43_1}, we have
\be
\label{p43_3}
|(g, \Delta v)| + \epsilon |  ( \Delta \zto , \Delta v )|
\le {\frac 12} \| \Delta v \|^2
+ \| g \|^2 + \epsilon  \|  \Delta \zto \|^2.
\ee
It follows from  \eqref{p43_1}-\eqref{p43_3} that, for all
$t \ge 0$,
$$
  {\frac d{dt}} \| \nabla v \|^2
 \le c \left (
 \| \nabla u \|^2 + \| u \|^p_p \right )
 + \epsilon c \left (
  \|   \Delta \zto  \|^2 +  \|   \Delta \zto  \|^p_p
\right ) + c
$$
\be
\label{p43_4}
 \le c \left (
 \| \nabla u \|^2 + \| u \|^p_p \right )
 + \epsilon c p_2(\theta_t \omega) + c,
\ee
where
$p_2(\theta_t \omega)=
  \|   \Delta \zto  \|^2 +  \|   \Delta \zto  \|^p_p
 $.
 Let $T(B, \omega)$ be the   constant in Lemma \ref{lem42},
 fix  $t \ge T(B,\omega)$ and $ s \in (t, t+1)$.
 Integrating  \eqref{p43_4} in
 $(s, t+1)$    we find that
 $$
 \| \nabla v(t+1, \omega, v_0(\omega)) \|^2
 \le \| \nabla v(s, \omega, v_0(\omega)) \|^2
 + \epsilon c \int_s^{t+1} p_2(\theta_\tau \omega) d \tau
 $$
 $$
 +
 c \int_s^{t+1}
 \left (
   \| \nabla u(\tau, \omega, u_0(\omega)) \|^2 + \|  u(\tau, \omega, u_0(\omega)) \|^p_p
 \right ) d \tau  +c.
 $$
 $$
  \le \| \nabla v(s, \omega, v_0(\omega)) \|^2
 + \epsilon c  \int_t^{t+1} p_2(\theta_\tau \omega) d \tau
 $$
 $$
 +
 c \int_t^{t+1}
 \left (
   \| \nabla u(\tau, \omega, u_0(\omega)) \|^2 + \|  u(\tau, \omega, u_0(\omega)) \|^p_p
 \right ) d \tau +c .
 $$
 Integrating the above with respect to $s$ in $(t, t+1)$, we have
 $$  \| \nabla v(t+1, \omega, v_0(\omega)) \|^2
  \le \int_t^{t+1}  \| \nabla v(s, \omega, v_0(\omega)) \|^2 ds
 + \epsilon c \int_t^{t+1} p_2(\theta_\tau \omega) d \tau
 $$
 $$
 +
 c \int_t^{t+1}
 \left (
   \| \nabla u(\tau, \omega, u_0(\omega)) \|^2 + \|  u(\tau, \omega, u_0(\omega)) \|^p_p
 \right ) d \tau +c.
 $$
 Now replacing $\omega$ by $\theta_{-t -1} \omega$, we get that
 $$  \| \nabla v(t+1, \theta_{-t -1} \omega, v_0(\theta_{-t -1} \omega)) \|^2
 $$
 $$
  \le \int_t^{t+1}
  \| \nabla v(s, \theta_{-t -1} \omega, v_0(\theta_{-t -1} \omega)) \|^2 ds
 + \epsilon c \int_t^{t+1} p_2(   \theta_{\tau -t -1} \omega) d \tau
 $$
\be
\label{p43_10}
 +
 c \int_t^{t+1}
 \left (
   \| \nabla u(\tau, \theta_{-t -1}\omega, u_0(\theta_{-t -1}\omega)) \|^2
+ \|  u(\tau, \theta_{-t -1}\omega, u_0(\theta_{-t -1}\omega)) \|^p_p
 \right ) d \tau +c .
\ee
Replacing $t$ by $t+1$ in Lemma \ref{lem41}, we find that there exists
$T_1=T_1(B, \omega)>0$, independent of $\epsilon$, such that for all $t \ge T_1$,
\be
\label{p43_11}
\int_t^{t+1} e^{\lambda (\tau -t-1)}
  \| \nabla v(\tau, \theta_{-t -1} \omega, v_0(\theta_{-t -1} \omega)) \|^2 d\tau
  \le c + \epsilon c r(\omega),
\ee
and
\be
\label{p43_12}
\int_t^{t+1} e^{\lambda (\tau -t-1)}
  \|   u (\tau, \theta_{-t -1} \omega, u_0(\theta_{-t -1} \omega)) \|^p_p d\tau
  \le c + \epsilon c r(\omega).
\ee
Since $e^{\lambda (\tau -t-1)} \ge e^{-\lambda}$ for $\tau \in (t, t+1)$, we obtain
from \eqref{p43_11}-\eqref{p43_12} that, for all
$t \ge T_1$,
\be
\label{p43_13}
\int_t^{t+1}  (
  \| \nabla v(\tau, \theta_{-t -1} \omega, v_0(\theta_{-t -1} \omega)) \|^2
   +
  \|   u (\tau, \theta_{-t -1} \omega, u_0(\theta_{-t -1} \omega)) \|^p_p
  ) d\tau
  \le ce^\lambda (1+   \epsilon   r(\omega)).
\ee
It follows from \eqref{p43_10}, \eqref{p43_13} and
Lemma  \ref{lem42}  that, there is $T_2 =T_2(B, \omega)>0$, independent of $\epsilon$,
such that for all $t\ge T_2$,
$$  \| \nabla v(t+1, \theta_{-t -1} \omega, v_0(\theta_{-t -1} \omega)) \|^2
 \le c_1 +  \epsilon c_2  r(\omega)
 + \epsilon c \int_{-1}^0 p_2(\theta_\tau \omega) d\tau
 $$
\be
\label{p43_14}
 \le c_1 + \epsilon  c_2  r(\omega)
 +  \epsilon c_3  \int_{-1}^0     e^{-{\frac \lambda 2}\tau} r(\omega)   d\tau
 \le c_1 + \epsilon c_4 r (\omega),
\ee
where we have used \eqref{zz}.
From \eqref{uandv} and \eqref{p43_14} we have, for all $t \ge T_2$,
\be
\label{p43_15}
 \| \nabla u  (t+1, \theta_{-t-1} \omega, u_0(\theta_{-t-1 } \omega)) \|^2
 \le c_5 + \epsilon c_6   r(\omega).
\ee
The lemma then follows from \eqref{p43_14} and \eqref{p43_15}.
\end{proof}

Next, we derive uniform estimates of solutions for large space and time variables.
Particularly, we show how these estimates depend on the small parameter $\epsilon$.

\begin{lem}
\label{lem44}
Let $0< \epsilon \le 1  $, $g \in L^2(\R^n)$ and
\eqref{f1}-\eqref{f4} hold.
Suppose $ B= \{B(\omega)\}_{\omega \in \Omega}\in \mathcal{D}$
and $u_0(\omega) \in B(\omega)$.
Then for every $\eta>0$ and  $P$-a.e. $\omega \in \Omega$,
there exist     $T = T  (B, \omega, \eta)>0$ and
$R = R  (\omega, \eta)>0$
such that    for all $t \ge T $,
$$
 \int_{|x| \ge R } | u(t, \theta_{-t } \omega, u_0(\theta_{-t  } \omega) )  (x) |^2
 dx \le \eta,
$$
where $T  (B, \omega, \eta)$ and $R(\omega, \eta)$ do not depend
on $\epsilon$.
\end{lem}

\begin{proof}
Let $\rho$ be a smooth function defined on $   \R^+$ such that
$0\le \rho(s) \le 1$ for all $s \in \R^+$, and
$$
\rho (s) = \left \{
\begin{array}{ll}
  0 & \quad \mbox{for} \ 0\le s \le 1; \\
 1 & \quad \mbox{for}  \  s \ge 2.
\end{array}
\right.
$$
Then there exists a positive  constant
$c$ such that
$ | \rho^\prime (s) | \le c$ for all $s \in \R^+$.

 Taking the inner product of \eqref{v1} with $\rho({\frac {|x|^2}{k^2}})v$
 in $\h$, we obtain that
 $$
 {\frac 12} {\frac d{dt}} \ii \rh |v|^2 dx
 +\lambda \ii \rh |v|^2 dx
 + \ii |\nabla v|^2 \rh   dx
 $$
 \be
 \label{p44_1}
=  \ii f(x, u) \rh v dx - \ii v \rhp {\frac {2x}{k^2}} \cdot \nabla v dx
 + \ii \left ( g + \epsilon \Delta \zto  \right ) \rh v dx.
 \ee
 By \eqref{f1} and \eqref{f2}, the first term on the right-hand side
 of \eqref{p44_1} satisfies
$$
  \ii f(x, u) \rh v dx
 =  \ii f(x, u) \rh u dx -  \epsilon  \ii f(x, u) \rh \zto dx
$$
$$
  \le
 -  {\frac 12} \alpha_1
\ii |u|^p \rh dx + \ii \psi_1 \rh dx
$$
\be
\label{p44_8}
+{\frac 12} \ii \psi_2^2 \rh dx
 + \epsilon c \ii  \left (  |\zto|^p +  | \zto|^2 \right ) \rh dx.
\ee
 Note that the second term on the right-hand side of
  \eqref{p44_1}  is bounded by
  $$
  | \ii v \rhp {\frac {2x}{k^2}} \cdot \nabla v dx|
    =
     |\int_{k\le |x| \le \sqrt{2} k} v \rhp {\frac {2x}{k^2}} \cdot \nabla v dx|
    $$
    \be
    \label{p44_8_a1}
     \le {\frac {2 \sqrt{2}}k} \int _{k \le |x| \le \sqrt{2} k} |v|
     \; |\rhp| \;  | \nabla v| dx
    \le {\frac ck} (\| v \|^2 + \| \nabla v \|^2 ).
  \ee
For the last term on the right-hand side of \eqref{p44_1}, we have
\be
\label{p44_9}
| \ii (g + \epsilon \Delta \zto) \rh v dx|
\le
{\frac 12} \lambda \ii \rh |v|^2 dx
+ {\frac 1\lambda} \ii (g^2 + \epsilon^2 |\Delta \zto |^2 ) \rh dx.
\ee
It follows from  \eqref{p44_1}-\eqref{p44_9}
that
$$
  {\frac d{dt}} \ii \rh |v|^2 dx
+   \lambda \ii \rh |v|^2 dx
 $$
$$
\le {\frac ck} (\| \nabla v \|^2 + \| v \|^2)
+
 c \ii
\left (   |\psi_1| +  |\psi_2|^2
+   g^2 \right )\rh
dx
$$
\be
\label{p44_11}
+ \epsilon c \ii \left (
|\Delta \zto |^2  + |\zto |^2 + | \zto |^p \right )
\rh dx.
\ee
Then using Lemmas \ref{lem41}-\ref{lem43} and following the process
of \cite{blw}, after detailed calculations we find that, given $\eta>0$,
 there exist $T=T(B, \omega, \eta)$ and $R=R(B, \eta)$,  which   are independent
 of $\epsilon$, such that   for all
$t\ge T $ and $k \ge R $,
$$
\int_{|x| \ge   k}   |v(t, \theta_{-t}\omega, v_0(\theta_{-t}\omega) )|^2 dx
\le
  \eta,
$$
which along with
\eqref{uandv} implies the lemma.
\end{proof}

\section{Upper semicontinuity of random attractors for Reaction-Diffusion equations
on $\R^n$}
\setcounter{equation}{0}

In this section, we prove the upper semicontinuity
of random attractors for the    Reaction-Diffusion equation defined on $\R^n$
when the stochastic perturbations approach zero.
To this end,  we first  establish the convergence of solutions
of problem \eqref{rd1}-\eqref{rd2} when $\epsilon \to 0$, and then
   show that the union
of all perturbed random attractors is precompact in $L^2(\R^n)$.

To indicate dependence of solutions on $\epsilon$,
in this section, we write the solution of problem \eqref{rd1}-\eqref{rd2}
as $u^\epsilon$, and the corresponding cocycle as $\Phi_\epsilon$.
Given $0< \epsilon \le 1$, it follows from Lemma \ref{lem42} that,
for every $B=\{B(\omega)\}_{\omega \in \Omega}  \in {\mathcal{D}}$
and $P$-a.e. $\omega \in \Omega$, there exists $T=T(B, \omega)>0$, independent of
$\epsilon$, such that for all $t \ge T$,
\be
\label{uniabs}
\| \Phi_\epsilon (t, \theta_{-t} \omega, B(\theta_{-t}\omega)) \|
\le M+ \epsilon M r(\omega),
\ee
where
$M$ is a positive deterministic  constant
independent of $\epsilon$, and $r(\omega)$ is the tempered
function in \eqref{z2}.
Denote by
\be
\label{rdeabs1}
K_\epsilon (\omega) = \{ u \in L^2(\R^n): \ \| u \| \le M+ \epsilon M r(\omega) \},
\ee
and
\be
\label{rdeabs2}
K  (\omega) = \{ u \in L^2(\R^n): \ \| u \| \le M+   M r(\omega) \},
\ee
where $M$ is the constant in \eqref{uniabs}.
Then  for every $0<\epsilon \le 1$, $\{K_\epsilon (\omega)\}_{\omega \in \Omega}$
 is a closed absorbing set
for $\Phi_\epsilon$ in $  {\mathcal{D}}$ and
\be
\label{sec6_1}
\bigcup\limits_{0< \epsilon \le 1} K_\epsilon (\omega)
\subseteq K(\omega).
\ee
It follows from the invariance of the random attractor
$\{{\mathcal{A}}_\epsilon (\omega) \}_{\omega \in \Omega}$ and
\eqref{sec6_1} that
\be
\label{sec6_2}
\bigcup\limits_{0< \epsilon \le 1} {\mathcal{A}}_\epsilon (\omega)
\subseteq \bigcup\limits_{0< \epsilon \le 1} K_\epsilon (\omega)
\subseteq K(\omega).
\ee
On the other hand,
by  Lemmas \ref{lem42} and \ref{lem43}, we find that,
for every  $0< \epsilon \le 1$
and $P$-a.e. $\omega \in \Omega$, there exists $T_1=T_1( \omega)>0$, independent of
$\epsilon$, such that for all $t \ge T_1$,
\be
\label{sec6_3}
\| \Phi_\epsilon (t, \theta_{-t} \omega, K(\theta_{-t}\omega)) \|_{H^1(\R^n)}
\le M_1 + \epsilon M_1 r(\omega)
\le M_1 + M_1r(\omega),
\ee
where $K(\omega)$ is given in \eqref{rdeabs2} and
$M_1$ is a positive deterministic  constant
independent of $\epsilon$.
By \eqref{sec6_2} and \eqref{sec6_3} we obtain that, for every $0<\epsilon \le 1$,
$P$-a.e. $\omega \in \Omega$ and $t \ge T_1$,
\be
\label{sec6_4}
\| \Phi_\epsilon (t, \theta_{-t} \omega, {\mathcal{A}}_\epsilon
(\theta_{-t}\omega)) \|_{H^1(\R^n)}
\le M_1 + M_1r(\omega).
\ee
By invariance,  ${\mathcal{A}}_\epsilon (\omega)
=\Phi_\epsilon (t, {\theta_{-t}} \omega, {\mathcal{A}}_\epsilon (\theta_{-t}
\omega))$ for all $t \ge 0$ and $P$-a.e. $\omega \in \Omega$.
Therefore, by \eqref{sec6_4} we have that,
for $P$-a.e. $\omega \in \Omega$,
\be
\label{sec6_5}
\| u \|_{H^1(\R^n)} \le M_1 + M_1 r(\omega),
\quad \forall \  u \in \bigcup\limits_{0< \epsilon \le 1}
 {\mathcal{A}}_\epsilon (\omega).
\ee
We remark that \eqref{sec6_5} is important for
proving the precompactness of the union
$\bigcup\limits_{0< \epsilon \le 1}
 {\mathcal{A}}_\epsilon (\omega)$
 in $L^2(\R^n)$.

\begin{lem}
\label{lem61}
Let  $g \in L^2(\R^n)$ and
\eqref{f1}-\eqref{f4} hold.
Then the union $\bigcup\limits_{0< \epsilon \le 1}
 {\mathcal{A}}_\epsilon (\omega)$
 is precompact in $L^2(\R^n)$.
\end{lem}

\begin{proof}
Given $\eta>0$, we want to show that the set
$\bigcup\limits_{0< \epsilon \le 1}
 {\mathcal{A}}_\epsilon (\omega)$
 has a finite covering of balls of radii less than
 $\eta$. Let $R$ be a positive number and denote by
 $$
 Q_R= \{x\in \R^n: \ |x| <R\}
 \quad \mbox{and} \quad  Q_R^c = \R^n \setminus Q_R.
 $$
 Let $\{K(\omega)\}_{\omega \in \Omega}$ be the random set
 given in \eqref{rdeabs2}.
 By Lemma \ref{lem44}, we find that,
 given  $\eta>0$ and  $P$-a.e. $\omega \in \Omega$,
there exist     $T = T  (\omega, \eta)>0$ and
$R = R  (\omega, \eta)>0$ (independent of $\epsilon$)
such that    for all $t \ge T $
and $u_0(\theta_{-t} \omega) \in K(\theta_{-t}\omega)$,
\be
\label{pl61_1}
 \int_{|x| \ge R } | u^\epsilon (t, \theta_{-t } \omega, u_0(\theta_{-t  } \omega) )  (x) |^2
 dx \le {\frac {\eta^2}{16}}.
\ee
 By \eqref{sec6_2}, $u_0(\theta_{-t} \omega) \in {\mathcal{A}}_\epsilon
   (\theta_{-t}\omega)$ implies
   that $u_0(\theta_{-t} \omega) \in K(\theta_{-t}\omega)$.
   Therefore it follows from
   \eqref{pl61_1} that, for every $0<\epsilon \le 1$,
   $P$-a.e. $\omega \in \Omega$,
   $t\ge T$  and $u_0(\theta_{-t} \omega) \in {\mathcal{A}}_\epsilon
   (\theta_{-t}\omega)$,
$$
 \int_{|x| \ge R } | u^\epsilon (t, \theta_{-t } \omega, u_0(\theta_{-t  } \omega) )  (x) |^2
 dx \le  {\frac {\eta^2}{16}},
$$
which along with
the invariance of $ \{{\mathcal{A}}_\epsilon (\omega)\}_{\omega \in \Omega}$
 shows that, for $P$-a.e. $\omega \in \Omega$,
 $$
 \int_{|x| \ge R } |  u  (x) |^2
 dx \le  {\frac {\eta^2}{16}},
  \quad \forall \ u \in \bigcup\limits_{0< \epsilon \le 1}
 {\mathcal{A}}_\epsilon (\omega),
$$
that is for $P$-a.e. $\omega$,
\be
\label{pl61_2}
\| u \|_{L^2(Q^c_R)} \le {\frac \eta{4}}, \quad
\forall \ u \in \bigcup\limits_{0< \epsilon \le 1}
 {\mathcal{A}}_\epsilon (\omega).
 \ee
 On the other hand, \eqref{sec6_5} implies that
 the set  $   \bigcup\limits_{0< \epsilon \le 1}
 {\mathcal{A}}_\epsilon (\omega)$
 is bounded in $H^1(Q_R)$ for $P$-a.e. $\omega \in \Omega$.
 By the compactness of embedding $H^1(Q_R)
 \subseteq  L^2(Q_R)$ we find that, for the given
 $\eta$, the set
 $   \bigcup\limits_{0< \epsilon \le 1}
 {\mathcal{A}}_\epsilon (\omega)$
 has a finite covering of balls of radii less than
 ${\frac \eta{4}}$ in $L^2(Q_R)$. This along with
 \eqref{pl61_2} shows that
 $   \bigcup\limits_{0< \epsilon \le 1}
 {\mathcal{A}}_\epsilon (\omega)$
 has a finite covering of balls of radii less than $\eta$
 in $L^2(\R^n)$.
\end{proof}

Next, we investigate the limiting behavior of solutions
of problem \eqref{rd1}-\eqref{rd2} when $\epsilon \to 0$.
We further assume that the nonlinear function
$f$ satisfies, for all $x \in \R^n$ and $s \in \R$,
\be
\label{fcond}
|{\frac {\partial f}{\partial s}} (x,s)|
\le \alpha_3 |s|^{p-2} + \psi_4 (x),
\ee
where $\alpha_3>0$,  $\psi_4 \in L^\infty(\R^n)$ if $p=2$, and
$\psi_4 \in L^{\frac p{p-2}} (\R^n)$ if $p>2$.

Under condition \eqref{fcond}, we  will show that, as $\epsilon \to 0$,
the solutions of the perturbed equation \eqref{rd1}
converge to  the limiting deterministic equation:
\be
\label{deter1}
{\frac {du}{dt}} + \lambda u -\Delta u =f(x,u) + g(x),
\quad x \in \R^n, \  t>0.
\ee

\begin{lem}
\label{lem62}
Suppose  $g \in L^2(\R^n)$,
\eqref{f1}-\eqref{f4} and \eqref{fcond} hold.
Given $0< \epsilon \le 1 $,  let $u^\epsilon$ and $u$ be the
solutions of equation \eqref{rd1} and \eqref{deter1}
with initial conditions $u^\epsilon_0$ and $u_0$, respectively.
Then for $P$-a.e. $\omega \in \Omega$ and $t\ge 0$, we have
$$\| u^\epsilon (t, \omega, u_0^\epsilon) - u(t, u_0)\|^2
\le
 c e^{ct}
\| u^\epsilon_0 - u_0 \|^2
+\epsilon c e^{c t}
\left ( r(\omega) + \|u^\epsilon_0 \|^2 + \|u_0 \|^2
\right ),
$$
where $c$ is a positive deterministic constant independent of $\epsilon$,
and $r(\omega)$ is the tempered function in \eqref{z2}.
\end{lem}

 \begin{proof}
 Let $v^\epsilon = u^\epsilon (t, \omega, u_0^\epsilon)-\epsilon z  (\theta_t \omega)$
 and $W= v^\epsilon -u$. Since $v$ and $u$ satisfy \eqref{v1}
 and \eqref{deter1}, respectively, we find that $W$ is a solution
 of the equation:
 $$
  {\frac {\partial W}{\partial t}}
  +\lambda W -\Delta W
  =f(x, u^\epsilon) - f(x, u) +\epsilon \Delta z(\theta_t \omega).
 $$
 Taking the inner product of the above with $W$ in $L^2(\R^n)$ we get
 \be
 \label{p62_1}
 {\frac 12} {\frac d{dt}} \| W\|^2
 +\lambda \| W\|^2 + \| \nabla W \|^2
 =\int_{\R^n} (f(x, u^\epsilon) - f(x, u)) W dx
 +\epsilon \int_{\R^n} \Delta z(\theta_t \omega) W dx.
 \ee
 For the first term on the right-hand side of \eqref{p62_1},
 by \eqref{f3} and \eqref{fcond} we have
 $$
 \int_{\R^n} (f(x, u^\epsilon) - f(x, u)) W dx
 =\int_{\R^n} {\frac {\partial f}{\partial s}} (x,s) (u^\epsilon -u) W dx
 $$
 $$
 = \int_{\R^n} {\frac {\partial f}{\partial s}} (x,s) W^2 dx
 +\epsilon \int_{\R^n} {\frac {\partial f}{\partial s}}(x,s) z (\theta_t \omega) W dx
 $$
 $$
 \le \beta \| W\|^2 +\epsilon \alpha_3 \int_{\R^n} (|u^\epsilon| + |u|)^{p-2}
 |z (\theta_t \omega)| |W| dx
 +\epsilon \int_{\R^n} \psi_4 |z (\theta_t \omega)| |W| dx
 $$
 \be
 \label{p62_2}
 \le \beta \| W\|^2 +\epsilon c
 \left ( \| u^\epsilon\|^p_p
 + \| u \|^p_p +
 \| z (\theta_t \omega) \|^p_p + \| W \|^p_p + \| \psi_4 \|^{\frac p{p-2}}_{\frac p{p-2}}
 \right ).
 \ee
 By the Young inequality, the last term on the right-hand side of
 \eqref{p62_1} is bounded by
 \be
 \label{p62_3}
 \epsilon \int_{\R^n} | \Delta z(\theta_t \omega) W | dx
 \le {\frac 12} \epsilon \| \Delta z(\theta_t \omega) \|^2
  + {\frac 12} \epsilon  \| W \|^2
  \le {\frac 12} \epsilon \| \Delta z(\theta_t \omega) \|^2
  + {\frac 12}   \| W \|^2.
  \ee
  It follows from \eqref{p62_1}-\eqref{p62_3} that
 $$ {\frac d{dt}} \| W \|^2
  \le c \| W \|^2
  + \epsilon c + \epsilon c
 \left ( \| u^\epsilon\|^p_p
 + \| u \|^p_p +
 \| z (\theta_t \omega) \|^p_p
 + \| \Delta z(\theta_t \omega)\|^2 + \| W \|^p_p
 \right )
 $$
 $$
  \le c \| W \|^2
  + \epsilon c + \epsilon c
 \left ( \| u^\epsilon\|^p_p
 + \| u \|^p_p +
 \| z (\theta_t \omega) \|^p_p
 + \| \Delta z(\theta_t \omega)\|^2
 \right )
$$
\be
 \label{p62_4}
  \le c \| W \|^2
  + \epsilon c + \epsilon c
 \left ( \| u^\epsilon\|^p_p
 + \| u \|^p_p  \right ) + \epsilon c
  e^{{\frac 12} \lambda |t|} r(\omega),
 \ee
 where we have used $W= u^\epsilon (t, \omega, u_0^\epsilon)-\epsilon z  (\theta_t \omega)-u$
 ,  the fact $0<\epsilon \le 1$ and \eqref{zz}.
 Integrating \eqref{p62_4} on (0,t) we obtain
 $$
 \| W(t) \|^2
 \le e^{c t} \| W(0) \|^2
 + \epsilon c
 + \epsilon c r(\omega)  e^{ct} \int_0^t e^{ ({\frac 12} \lambda     -c)s} ds
 $$
 $$
 +\epsilon c  \int_0^t e^{c(t-s)}
 \left ( \| u^\epsilon (s, \omega, u^\epsilon_0) \|^p_p
 + \| u(s,u_0) \|^p_p \right ) ds
 $$
\be
\label{p62_4a1}
 \le e^{c t} \| W(0) \|^2
 + \epsilon c_1
 + \epsilon c_1 r(\omega)  e^{c_2 t}
 +\epsilon c e^{ct}  \int_0^t
 \left ( \| u^\epsilon (s, \omega, u^\epsilon_0) \|^p_p
 + \| u(s,u_0) \|^p_p \right ) ds.
 \ee
 It follows from \eqref{p41_8} that
 $$
   \int_0^t e^{\lambda (s -t)} \| u^\epsilon (s, \omega, u_0^\epsilon)\|^p_p ds
\le e^{-\lambda t} \| v^\epsilon _0 (\omega) \|^2
+ \epsilon c  \int_0^t e^{\lambda(s -t)}  p_1(\theta_s \omega ) ds
+ c,
$$
which together with \eqref{p41_6a1}
implies that, for all $ t \ge 0$,
$$
   \int_0^t e^{\lambda  s  } \| u^\epsilon (s, \omega, u_0^\epsilon)\|^p_p ds
\le   \| v_0^\epsilon (\omega) \|^2
+ \epsilon c  \int_0^t e^{\lambda s  }  p_1(\theta_s \omega ) ds
+ c e^{\lambda t}
$$
\be
\label{p62_5}
\le   \| v_0^\epsilon (\omega) \|^2
+   c r(\omega) \int_0^t e^{{\frac 32} \lambda s  }    ds
+ c e^{\lambda t}
\le    \| u_0^\epsilon -\epsilon z(\omega) \|^2
+   c_3 r(\omega)   e^{c_4 t  }
+ c e^{\lambda t}.
\ee
Since $e^{\lambda s} \ge 1$ for all $s \in [0,t]$, we obtain from
\eqref{p62_5} that
\be
\label{p62_6}
   \int_0^t   \| u^\epsilon (s, \omega, u_0^\epsilon)\|^p_p ds
\le   2 \| u_0^\epsilon   \|^2 + 2  \| z(\omega) \|^2
+   c_3 r(\omega)   e^{c_4 t  }
+ c e^{\lambda t}.
\ee
Similarly, by \eqref{deter1} for $\epsilon =0$, we can also get that
\be
\label{p62_7}
   \int_0^t   \| u  (s,   u_0 )\|^p_p ds
\le   c \| u_0    \|^2
+ c e^{\lambda t}.
\ee
By \eqref{z2}, \eqref{p62_4a1} and \eqref{p62_6}-\eqref{p62_7} we find that,
for all $t \ge 0$,
\be
\label{p62_10}
\| W(t) \|^2
\le e^{ct} \| W(0) \|^2
+\epsilon c e^{c_5t}
\left ( r(\omega) + \|u^\epsilon_0 \|^2 + \|u_0 \|^2
\right ).
\ee
Finally, by \eqref{zz} and  \eqref{p62_10} we have, for all $t \ge 0$,
$$\| u^\epsilon (t, \omega, u_0^\epsilon) - u(t, u_0)\|^2
=\| W(t) + \epsilon z(\theta_t \omega) \|^2
\le 2 \|W(t)\|^2 + c_6 \epsilon e^{c_7t} r(\omega)
$$
$$
\le 2 e^{ct} \|W(0) \|^2
+\epsilon c e^{c_8t}
\left ( r(\omega) + \|u^\epsilon_0 \|^2 + \|u_0 \|^2
\right )
$$
$$
\le
 2 e^{ct}
\| u^\epsilon_0 - u_0 - \epsilon z(\omega) \|^2
+\epsilon c e^{c_8t}
\left ( r(\omega) + \|u^\epsilon_0 \|^2 + \|u_0 \|^2
\right )
$$
$$
\le
 4 e^{ct}
\| u^\epsilon_0 - u_0 \|^2
+\epsilon c_9 e^{c_8t}
\left ( r(\omega) + \|u^\epsilon_0 \|^2 + \|u_0 \|^2
\right ).
$$
This completes the proof.
 \end{proof}

 We are now in  a position to establish
 the upper semicontinuity of the perturbed
 random attractors for problem \eqref{rd1}-\eqref{rd2}.

 \begin{thm}
Let $g \in L^2(\R^n)$,
\eqref{f1}-\eqref{f4} and \eqref{fcond} hold.
Then for $P$-a.e. $\omega \in \Omega$,
\be
\label{thm63}
\lim_{\epsilon \to 0} dist_{L^2(\R^n)}
({\mathcal{A}}_\epsilon (\omega), {\mathcal{A}})
=0,
\ee
where
$$dist_{L^2(\R^n)} ({\mathcal{A}}_\epsilon (\omega), {\mathcal{A}})
= \sup_{a \in {\mathcal{A}}_\epsilon (\omega)} \inf_{b \in {\mathcal{A}}}
\| a -b \|_{L^2(\R^n)}.
$$
\end{thm}

\begin{proof}
Note that $\{K_\epsilon (\omega)\}_{\omega \in \Omega}$ is a closed
absorbing set
for $\Phi_\epsilon$ in ${\mathcal{D}}$, where
$K_\epsilon (\omega)$ is given by \eqref{rdeabs1}.
By \eqref{rdeabs1} we find that
\be
\label{p63_1}
\limsup_{\epsilon \to 0} \| K_\epsilon (\epsilon) \| \le M,
\ee
where $M$ is the positive deterministic constant in \eqref{rdeabs1}.
Let $\epsilon_n \to 0$ and $u_{0,n} \to u_0$ in $L^2(\R^n)$.
Then by Lemma \ref{lem62} we find that, for $P$-a.e. $\omega \in \Omega$
and $t \ge 0$,
\be
\label{p63_2}
\Phi_{\epsilon_n} (t, \omega, u_{0,n}) \to \Phi (t, u_0).
\ee
Notice that \eqref{p63_1}-\eqref{p63_2} and Lemma \ref{lem61}
indicate   all
  conditions \eqref{upp1}-\eqref{upp3}
  are satisfied, and hence \eqref{thm63}
follows  from Theorem \ref{semicont} immediately.
\end{proof}

 \end{document}